\def\th{{^{\rm th}}}
 \def\D{{\partial}}
 \def\BCH{{\rm BCH}}
 
 \def\Aut{{\rm Aut}}
 \def\Diff{{\rm Diff}}
 \def\ad{{\rm ad}}

 \def\takes{{\colon}}
 \def\QED{{\hfill$\Box$}}
 \def\rats{{\bf Q}}
\documentclass{gtart}
\usepackage{enumerate}
\usepackage{latexsym}
\usepackage{amsfonts}
\usepackage{amssymb}
\usepackage{amscd}
\usepackage{amsbsy}
\usepackage{amsmath}
\usepackage{tikz}
\usetikzlibrary{arrows.meta}

\begin{document}

\title{An explicit symmetric DGLA model of a bi-gon}

\authors{Nir Gadish, Itay Griniasty and Ruth Lawrence}

\address{University of Chicago, Weizmann Institute and Hebrew University}
\email{nirg@math.uchicago.edu, itay.griniasty@weizmann.ac.il, ruthel@ma.huji.ac.il}

\begin{abstract}

We give explicit formulae for a DGLA model of the bi-gon which is symmetric under the geometric symmetries of the cell. This follows the work of Lawrence-Sullivan on the (unique) DGLA model of the interval and its construction uses deeper knowledge of the structure of such models and their localisations for non-simply connected spaces.

 \end{abstract}
 \primaryclass{17B55}\secondaryclass{17B01, 55U15}
\keywords{DGLA, infinity structure, Maurer-Cartan, Baker-Campbell-Hausdorff formula}
 \maketitlepage

 \section {Introduction}
For a regular cell complex $X$, it is possible to associate a DGLA model $A=A(X)$ satisfying the following conditions
 \begin{itemize}
 \item as a Lie algebra, $A(X)$ is freely generated by a set of generators, one for each cell in $X$ and whose grading is one less than the geometric degree of the cell;
 \item vertices (that is $0$-cells) in $X$ give rise to generators $a$ which satisfy the Maurer-Cartan equation $da+{1\over2}[a,a]=0$ (a flatness condition);
 \item for a cell $x$ in $X$, the part of $\D{x}$ without Lie brackets is the geometric boundary $\D_0x$ (where an orientation must be fixed on each cell);
 \item (locality) for a cell $x$ in $X$, $\D{x}$ lies in the Lie algebra generated by the generators of $A(X)$ associated with cells of the closure $\bar{x}$.
 \end{itemize}

 The existence and general construction of such a model was demonstrated by Sullivan in the Appendix to [4], and its uniqueness up to DGLA isomorphism (along with other general properties) were proved in [1].

 In this paper we give an explicit model for the bi-gon which exhibits the dihedral symmetry of the bi-gon. The existence of such a symmetric model, though known to be not unique, the first four orders and various other properties of such a model were given in [2], the current manuscript is the first time that an explicit closed formula has been given for one. In this section we collect some general facts about models of cell complexes (see [3]) while in succeeding sections we focus on the bi-gon and its boundary.

{\bf Points and localisation.} An element $a\in{}A_{-1}$ is called a {\it point} in the model, if it satisfies the Maurer-Cartan equation $\D{a}+{1\over2}[a,a]=0$. From any point in the model, we can define a {\it twisted differential} $\D_a$ by $\D_a=\D+\ad_a$. Then $\D_a$ is a derivation (immediate from Leibnitz and Jacobi) while
$$\D_a^2=\D^2+[\D,\ad_a]+\ad_a^2=\D^2+\ad_{\D{a}}+\frac{1}{2}[\ad_a,\ad_a]
=\D^2+\ad_{\D{a}+\frac{1}{2}[a,a]}=0\>.$$
By the {\sl  localisation} of $A$ to $a$, denoted $A(a)$, we will mean the complex
$$0\longleftarrow\ker\D_a|_{A_0}\longleftarrow{}A_1\longleftarrow{}A_2\longleftarrow\cdots$$
under the action of the differential $\D_a$.

{\bf Edges and flows.} Any element $e\in{}A_0$ defines a {\it flow} on $A$ by
 $$\frac{dx}{dt}=\D{e}-\ad_e(x)\quad{\rm on}\quad A_{-1}\>,\qquad
 \frac{dx}{dt}=-\ad_e(x)\quad{\rm on}\quad A_{\geq0}\>,$$
This flow is called the flow {\it by} $e$, and preserves the grading. In grading $-1$, it preserves flatness, meaning that if the initial condition is at a point ($x(0)$ satisfies Maurer-Cartan) then at all time its value also satisfies Maurer-Cartan.  Linearity of the differential equation in $e$ ensures that flowing by $e$ for time $t$ is equivalent to flowing by $te$ for a unit time. Denote the result of flowing by $e$ from $a$ for unit time, by $u_e(a)$, so that the solution of the above differential equation satisfies $x(t)=u_{te}(x(0))$.

\noindent{\bf Definition}$\>${\sl Let $\BCH(x,y)$ (given by the Baker-Campbell-Hausdorff formula) denote the unique element in the free Lie algebra on the two generators $x$ and $y$ for which  $(\exp{x}).(\exp{y})=\exp{\BCH(x,y)}$ in the universal enveloping algebra of $A$, or equivalently $$(\exp{\ad_x})\circ(\exp{\ad_y})=\exp{\ad_{\BCH(x,y)}}\in\Aut(A)\>.$$}

\noindent{\bf Properties}
 \begin{itemize}
 \item $\BCH(x,y)$  is a sum of terms starting with $x+y$ with the remaining terms all being Lie brackets of increasing order of $x$'s and $y$'s,
     \begin{align*}
    \BCH(x,y)=x+y&+{1\over2}[x,y]+{1\over12}(X^2y+Y^2x)-{1\over24}XYXy\\
    -{1\over720}(X^4y+Y^4x)&+{1\over120}(X^2Y^2x+Y^2X^2y)+{1\over360}(XY^3x+YX^3y)+\cdots
     \end{align*}
     where $X,Y$ denote $\ad_x,\ad_y$.  In particular for commuting $x,y$, $\BCH(x,y)=x+y$ so that for example $\BCH(x,x)=2x$.
 \item Uniqueness implies associativity of BCH, that is, $\BCH\big(\BCH(x,y),z\big)=\BCH\big(x,\BCH(y,z)\big)$ for any symbols $x,y,z$. Denote the combined BCH of $n$ symbols by $\BCH(x_1,\ldots,x_n)$ so that $$\exp{\BCH(x_1,\ldots,x_n)}=(\exp{x_1})\cdots(\exp{x_n})\>,$$
     in $U(A)$ and $$\exp{\BCH(\ad_{x_1},\ldots,\ad_{x_n})}
     =(\exp{\ad_{x_1}})\cdots(\exp{\ad_{x_n}})\in\Aut(A)\>.$$
 \item Uniqueness similarly implies that $\BCH(x,-x)=0$ while $BCH(-x,-y)=-\BCH(y,x)$.
 \item The Jacobi identity implies that $\BCH(\ad_x,\ad_y)=\ad_{\BCH(x,y)}$.
\item $\BCH(x,y,-x)=(\exp{\ad_x})y$
\item $\BCH(\exp(\ad_e)x,\exp(\ad_e)y))=\exp(\ad_e)\BCH(x,y)$.
\end{itemize}

A flow by $e$ for unit time followed by a flow by $f$ for unit time is equivalent to a flow by $\BCH(e,f)$ for unit time,
 $$u_f\circ{}u_e=u_{\BCH(e,f)}\>,$$
at all gradings [3], so that $e\mapsto{}u_e$ is a homomorphism $(A^{(0)},\BCH)\longrightarrow\Diff(A)$.

\noindent{\bf Lemma 1}$\>${\sl For any flow $e$ from $a$ to $b$, $\D_b\circ\exp(-\ad_e)=\exp(-\ad_e)\circ\D_a$. That is, $\exp(-\ad_e)$ intertwines the localisation $A(a)$ to the localisation $A(b)$.}

\noindent{\bf Proof}$\>$ For arbitrary $n\geq0$, let $x(t)\in{}A_{-1}$ and $v(t)\in{}A_n$ be flows by $e$ starting with initial conditions $x(0)=a$ and $v(0)=u$. That is, assume $\dot{x}=\D{e}-\ad_e(x)$ and $\dot{v}=-\ad_e(v)$. Then we compute
 \begin{align*}
  \frac{d}{dt}\left(\D_{x(t)}v(t)\right)
  &=\frac{d}{dt}\left(\D{v(t)}+[x(t),v(t)]\right)\\
  &=\D{\dot{v}}+[\dot{x},v]-[x,\dot{v}]\\
  &=-\D(\ad_ev)+[\D{e}-\ad_ex,v]+[x,\ad_ev]\\
  &=-\left([\D{e},v]+[e,\D{v}]\right)+[\D{e},v]-\ad_e[x,v]\\
  &=-\ad_e\left(\D{v}+[x,v]\right)=-\ad_e\left(\D_x{v}\right)
 \end{align*}
Thus $\D_xv\in{}A_{n-1}$ also infinitesimally transforms by $-\ad_e$, that is it also flows by $e$.  Since $e$ is constant, the solution to such a flow is $\exp(-t.\ad_e)$ times the initial condition. Thus $v(1)=\exp(-\ad_e)v(0)$ while $\D_{x(1)}v(1)=\exp(-\ad_e)\D_{x(0)}v(0)$. Since this holds for all $v(0)\in{}A_n$, thus, recalling that $x(0)=a$ and $x(1)=b$, we have the equality of operators $A_n\longrightarrow{}A_{n-1}$ given by
$$\D_b\circ\exp(-\ad_e)=\exp(-\ad_e)\circ\D_a\>.$$\QED

\noindent{\bf Definition}$\>${\sl For any path $\gamma$, define $\BCH(\gamma)\in{}A_0$ by  $\BCH(\gamma)=\BCH(e_1,\ldots,e_m)$ where the path $\gamma$ consists of the edges $e_1,\ldots,e_m$ in that order.}

By [3], the differential on 1-cells is determined by the condition that a flow for unit time by the generator corresponding to the 1-cell takes the starting point to the ending point of the 1-cell. Explicitly
$$\D{e}=(\ad_e)b+\sum_{i=0}^\infty{\frac{B_i}{i!}}(\ad_e)^i(b-a)={E\over1-e^E}a+{E\over1-e^{-E}}b\>,$$
where $B_i$ denotes the $i\th$ Bernoulli number defined as
coefficients in the expansion
${x\over{e^x-1}}=\sum\limits_{n=0}^\infty{}B_n {x^n\over{}n!}$, $E\equiv\ad_e$ and the expressions in $E$ are considered as formal power series.

By [1], the set of equivalence classes of points in $A(X)$ under the equivalence relation generated by flow, is $\pi_0(X\coprod{*})$. More precisely, if $X$ has $c$ connected components and $\{a_1,\ldots,a_c\}$ is a choice of (base-)points, one in each connected component, then  the set of points in $A(X)$ is
$$\bigcup_{i=1}^c\big\{u_e(a_i)\bigm|e\in{}A_0\big\}\cup
\big\{u_e(0)\bigm|e\in{}A_0\big\}\>.$$
The last of these pieces (the connected component of $0$) is a copy of $A_0$, in the sense that the map $\pi_0\takes{}e\mapsto{}u_e(0)$ is injective. Meanwhile, for each $i$, the map $\pi_i\takes{}e\mapsto{}u_e(a_i)$ is a `fibration', in the sense that
 $$\BCH(\cdot,e)\takes\pi_i^{-1}(a_i)\longrightarrow\pi_i^{-1}(u_e(a_i))$$
describes explicitly a bijection between the fibres over $a_i$ and a general point $u_e(a_i)$, while the fibre over the specific point $a_i$ is generated as a vector space by $$\pi_i^{-1}(a_i)=\{\BCH(\gamma)|\gamma\in\pi_1(X,a_i)\}\>,$$
or equivalently as a Lie algebra by $\BCH(\gamma)$ where $\gamma$ ranges over a set of generators of $\pi_1(X,a_i)$. Thus the $i\th$ connected component of the set of points of $A(X)$ can be considered as a quotient of $A_0$ by the fundamental group of the associated component of $X$.

\section{The bi-gon}

Let $\bar{X}^2$ be the bi-gon with two vertices, two edge and one 2-cell. This has a model $\bar{A}^2$ with five generators as a free Lie algebra, $a,b$ (vertices) of grading $-1$, $e,f$ (1-cells) of grading 0 and $g$ (2-cell) of grading 1. Suppose that the orientations are such that $\D_0e=b-a$, $\D_0f=a-b$ and $\D_0g=e+f$. Let $X^2$ be $\bar{X}^2$ with the 2-cell removed, that is a circle with two vertices. Note that just as the model of an interval is unique, so is the model of any 1-complex, in particular of $X^2$.

The symmetries of $\bar{X}^2$ are given by the Klein four-group $V_4$ with four elements. This can be generated by two generators, the rotation $\sigma\takes{}a\leftrightarrow{}b$, $e\leftrightarrow{}f$ leaving $g$ fixed, and the reflection $\iota\takes{}e\leftrightarrow-f$ with $g\mapsto-g$ and leaving $a,b$ fixed.

The 2-cell with one vertex $\bar{X}^1$, has a model $\bar{A}^1$ with one generator in each degree $-1$,0,1, say $a$, $e$, $g$ respectively with $\D_0e=0$, $\D_0g=e$. The explicit model is
$$\D{}e=[e,a]\>,\qquad\D{}g=e-[a,g]\>.$$
Equivalently, $\D_ae=0$ and $\D_ag=e$. The only symmetry of $\bar{X}^1$ is given by inverting the orientation, $\iota\takes{}e\mapsto-e$, $g\mapsto-g$ and clearly the model is invariant under this symmetry. Note that $A^1=A(X^1)$ is generated by $a$,$e$ while its set of points is
$$\big\{u_{te}(a)\bigm|t\in{}\rats\big\}\cup
\big\{u_{te}(0)\bigm|t\in\rats\big\}=\{a\}\cup\big\{u_{te}(0)\bigm|t\in\rats\big\}\>.$$

Using the functoriality of the construction $X\mapsto{}A(X)$ under subdivision of intervals, one obtains a model of $\bar{X}^2$ (Figure 1). The action of $\D$ on the generators $a,b,e,f$ is fixed by the previous discussion as in $A^2$, so it only remains to give $\D{g}$. The model induced by subdividing $e$ into two intervals $e,f$ has
 $$\D{}g=\BCH(e,f)-[a,g]\eqno{(1)}$$
This is invariant under $\iota$ since $\BCH(-f,-e)=-\BCH(e,f)$ while $\iota(a)=a$. However it is not invariant under $\sigma$, under which it transforms to a model with
$$\D{}g=\BCH(f,e)-[b,g]\eqno{(2)}$$ We could describe these two models (1) and (2), as `based' at $a$ and $b$, respectively. The aim of this work is to provide an explicit model of $\bar{X}^2$ which is symmetric under the full dihedral action, this will be done by producing a model `based' at a symmetric point $x$.

\section{A symmetric model of the bi-gon}

The circle with two points $X^2$ has a (unique) model $<a,b,e,f>$ with the differential  determined (by Leibnitz) by its value on generators, which by definition on $a$, $b$ are
$$\D{a}=-{1\over2}[a,a]\>,\qquad{}\D{b}=-{1\over2}[b,b]\>,$$
while on edges $e$, $f$ are (as in the interval [3])
\begin{align*}
\D{e}&={E\over1-e^E}a+{E\over1-e^{-E}}b\>,\\
\D{f}&={F\over1-e^F}b+{F\over1-e^{-F}}a\>,
\end{align*}
where $E=\ad_e$, $F=\ad_f$. These ensure that $u_e(a)=b$ and $u_f(b)=a$. In particular $u_{\BCH(e,f)}(a)=a$ and indeed $\{w|u_w(a)=a\}=\ker\D_a$ is the one-dimensional vector space generated by $\BCH(e,f)$.

Thus $e$ and $-f$ are two different `directions' which both flow from $a$ to $b$ in unit time. More generally the condition that $h\in{}A_0^2$ flows from $a$ to $b$ in unit time is
\begin{align*}
u_h(a)=b&\Leftrightarrow{}u_h(a)=u_e(a)\\
&\Leftrightarrow{}u_{-e}\circ{}u_h(a)=a\\
&\Leftrightarrow{}u_{\BCH(h,-e)}(a)=a\\
&\Leftrightarrow{}\BCH(h,-e)\in\ker\D_a\\
&\Leftrightarrow{}\exists{}t\in{\bf R}:\BCH(h,-e)=t\BCH(e,f)\\
&\Leftrightarrow{}\exists{}t\in{\bf R}:h=\BCH(t\BCH(e,f),e)
\end{align*}
\noindent Here $t=0$ gives $h=e$ while $t=-1$ gives $h=-f$. To generate a symmetric `centre point' in $A_{-1}^2$ we travel from $a$ for time $1\over2$ in the direction midway between $e$ and $-f$ as supplied above by $h$ taking $t=-{1\over2}$.

\begin{figure}[h!]
	\centering
	\begin{tikzpicture}
	\begin{scope}[>=Latex]
	\draw (0,0) arc(180:-60:0.6cm);
	\draw (0,0) arc(-180:-60:0.6cm)[->];	
	\draw (0,0) node[anchor=east]{$a$};
	\filldraw (0,0) circle (1.2pt);
	\draw (0.6,0) node {\(g\)};
	\draw (0.6,0.6) node[anchor=south] {\(e\)};
	\draw (0,-0.85) node {\(\)};
	\end{scope}
	\end{tikzpicture}
	\qquad
	\begin{tikzpicture}
	\begin{scope}[>=Latex]
	\draw (0,0) arc (150:90:1.1cm) ;
	\draw (1.1*0.866025,0.55) arc (90:30:1.1cm)[<-];
	\draw (0.9,0.55) node[anchor=south]{$ f $};
	\draw (0,0) arc (-150:-90:1.1cm)[->];
	\draw (1.1*0.866025,-0.55) arc (-90:-30:1.1cm);	
	\draw (0.9,-0.55) node[anchor=north]{$ e $};		
	\filldraw (0,0) circle (1.2pt);
	\draw (0,0) node[anchor=east]{$ a $};
	\draw (1.9,0) node[anchor=west]{$ b $};
	\filldraw (1.1*2*0.866025,0) circle (1.2pt);	
	\draw (1.1*0.866025,0) node {$g$};
	\end{scope}
	\end{tikzpicture}
	\qquad
	\begin{tikzpicture}
	\begin{scope}[gray!50]
	\draw (0,0) arc (90+11:90-11:5cm);
	\draw (0,0) arc (115:65:2.25cm);
	\draw (0,0) arc (90+40:90-40:1.48cm);
	\draw (0,0) arc (90+52:90-52:1.2cm);
	
	\draw (0,0) arc (-90-11:-90+11:5cm);
	\draw (0,0) arc (-115:-65:2.25cm);
	\draw (0,0) arc (-130:-50:1.48cm);
	\draw (0,0) arc (-142:-38:1.2cm);
	\end{scope}
	
	\begin{scope}[>=Latex]
	\draw (0,0) arc (150:90:1.1cm) ;
	\draw (1.1*0.866025,0.55) arc (90:30:1.1cm)[<-];
	\draw (0,0) arc (-150:-90:1.1cm)[->];
	\draw (1.1*0.866025,-0.55) arc (-90:-30:1.1cm);	
	\filldraw (0,0) circle (1.2pt);
	\filldraw (1.1*2*0.866025,0) circle (1.2pt);	
	\draw [dashed] (1.1*0.866025,0) -- (1.1*2*0.866025,0) ;
	\draw (0,0) -- (1.1*0.8,0)[->] ;	
	\filldraw (1.1*0.866025,0) circle (1.2pt);		
	\draw (1.1*0.866025,0) node[anchor=north]{$x$};	
	\draw (0.7,-0.1) node[anchor=south]{\(v/2\)};
	\draw (0.9,-0.7) node[anchor=north]{$  $};	
	
	\draw (0,0) node[anchor=east]{$ a $};
	\draw (1.9,0) node[anchor=west]{$ b $};
	\draw (0.9,-0.55) node[anchor=north]{$ e $};
	\draw (0.9,0.55) node[anchor=south]{$ f $};	
	\end{scope}
	\end{tikzpicture}
	
	\caption{ {\it Left,}
		The complex \(\bar{X}^1\), can be sub-divided into  \(\bar{X}^2\) ({\it center}), however the derived algebra would not be symmetric under \(\sigma\).
		{\it Right,}
		the symmetric model \(\bar{A}^2\) is based at the point \(x\), half way between \(a\) and \(b\) along a flow according to \(v\) which interpolates symmetrically \(e\) and \(-f\). The gray lines describe flows along asymmetric interpolations of \(e\) and \(-f\).	}
\end{figure}
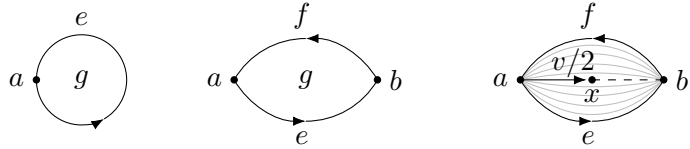

\noindent{\bf Lemma 2} Let $v=\BCH(-{1\over2}\BCH(e,f),e)$. Then $x=u_{v\over2}(a)\in{}A_{-1}^2$
 is a symmetric point in the model of $X_2$.

\noindent{\bf Proof}$\>$ That $x$ is a point follows from its form as the endpoint of a flow from $a$ (a point). Note that
$$u_v(a)=u_e\circ{}u_{-{1\over2}\BCH(e,f)}(a)=u_e(a)=b\>,$$
since $-{1\over2}\BCH(e,f)\in\ker\D_a$. To check symmetry, observe that
\begin{align*}
\iota(v)&=\BCH(-{1\over2}\BCH(-f,-e),-f)\\
&=\BCH({1\over2}\BCH(e,f),-f)\\
&=\BCH(-{1\over2}\BCH(e,f),\BCH(e,f),-f)\\
&=\BCH(-{1\over2}\BCH(e,f),e)=v
\end{align*}
so that $x=u_{v\over2}(a)$ is also preserved since $\iota(a)=a$. Under $\sigma$, $a\leftrightarrow{}b$ and $e\leftrightarrow{}f$ so that
\begin{align*}
\sigma(v)&=\BCH(-{1\over2}\BCH(f,e),f)\\
&=\BCH(-{1\over2}\exp(\ad_f)\BCH(e,f),\exp(\ad_f)f)\\
&=\exp(\ad_f)\BCH(-{1\over2}\BCH(e,f),f)\\
&=\BCH(f,-{1\over2}\BCH(e,f))=-\iota(v)=-v
\end{align*}
where in the second line we used that $f=\exp(\ad_f)f$. From this it follows that
$$\sigma(x)=u_{\sigma(v)\over2}(\sigma(a))=u_{-{v\over2}}(b)=u_{v\over2}(a)=x\>,$$
since $u_v(a)=b$.\QED

To obtain a model of $\bar{X}^2$ it suffices to add to the model of $X^2$, the additional generator $g$ along with its differential as
$$\D{g}=q-[x,g]\eqno{(3)}$$
where $x$ is the symmetric point defined above and $q\in{}A_0^2$ has zeroth order part (without Lie brackets) $\D_0g=e+f$. The requirement that $\D^2=0$ is satisfied so long as it is satisfied on generators. Since already $\D^2=0$ on $A^2$, it suffices to ensure that $\D^2g=0$. Now
$$\D^2g=\D(q-[x,g])=\D{}q-[\D{}x,g]+[x,\D{}g]=\D{}q-[\D{}x,g]+[x,q]-[x,[x,g]]=\D_xq\>,$$
since $[\D{}x,g]+[x,[x,g]]=[\D{}x+{1\over2}[x,x],g]=0$ as $x$ is a point. Thus (3) extends $A^2$ to a model of $\bar{A}^2$ on condition that $q$ is an element of grading 0 localised at the point $x$, that is, a solution of $\D_xq=0$, whose zeroth order (without Lie brackets) part is $e+f$. To obtain a symmetric model we additionally require that $q$ possesses suitable symmetry, $\iota(q)=-q$ and $\sigma(q)=q$ (recall that $\iota$ inverts the orientation of $g$ while $\sigma$ preserves it).

By Lemma~1, there is a bijection $\ker\D_a\longrightarrow\ker\D_{u_w(a)}$ given by
$$y\mapsto\BCH(-w,y,w)=\exp(-\ad_w)y\>,$$
while $\ker\D_a=\rats\BCH(e,f)$. Thus
$$\ker\D_x=\ker\D_{u_{v\over2}(a)}=\rats\BCH(-{v\over2},e,f,{v\over2})\>.$$
The zeroth order part of $\BCH$ is just the sum and thus $q=\BCH(-{v\over2},e,f,{v\over2})$ is localised at $x$ with zeroth order part $e+f$. Finally it is necessary to verify symmetry. Recall that $\iota(v)=v$ while $\sigma(v)=-v$. So
\begin{align*}
\iota(q)&=\BCH(-{v\over2},-f,-e,{v\over2})\\
&=-\BCH(-{v\over2},e,f,{v\over2})=-q\>,\\
\sigma(q)&=\BCH({v\over2},f,e,-{v\over2})\\
&=\BCH(-{v\over2},v,f,e,-v,{v\over2})\\
&=\BCH(-{v\over2},e,f,{v\over2})=q
\end{align*}
where in the last line we used that $\BCH(v,f,e,-v)=\BCH(e,f)$. This is because $v=\BCH({1\over2}\BCH(e,f),-f)$ and $-v=\BCH(-e,{1\over2}\BCH(e,f))$ and so
$$\BCH(v,f,e,-v)=BCH({1\over2}\BCH(e,f),-f,f,e,-e,{1\over2}\BCH(e,f))=\BCH(e,f)\>.$$

\noindent{\bf Theorem} The unique model $A^2=<a,b,e,f>$ of $X^2$ can be extended to a symmetric model $\bar{A}^2=<a,b,e,f,g>$ of $\bar{X}^2$ in which
$$\D{}g=q-[x,g]\>,$$
where $v=\BCH(-{1\over2}\BCH(e,f),e)$, $q=\BCH(-{v\over2},e,f,{v\over2})$ and $x=u_{v\over2}(a)$ are elements of $A^2$.

This model can be considered as `based' at the symmetric point $x\in{}A^2$ in the sense that $\D_xg=q\in<a,b,e,f>$

\section{Computations}

The first few orders of this explicit symmetric model can be computed using the formula for $\BCH$. We obtain
\begin{align*}
v&=\BCH(-{1\over2}\BCH(e,f),e)\\
&={1\over2}(e-f)+{1\over48}(E^2f-F^2e)+\cdots
\end{align*}
where odd orders (in the number of Lie brackets) vanish. From this,
$$\D{}v=b-a+{1\over4}(E-F)(a+b)+{1\over48}(E-F)^2(b-a)+{1\over96}[E+F,[E,F]](a+b)+\cdots$$
and then one obtains the central point in $A^2$ as
\begin{align*}
x&=u_{v\over2}(a)={1\over2}(a+b)+{1\over16}(F-E)(a-b)\\
&+{1\over3072}(E-F)^3(a-b)+{1\over384}[E+F,[E,F]](b-a)+\cdots
\end{align*}
where the term involving two Lie brackets vanishes. Finally we obtain
\begin{align*}
q&=\BCH(-{v\over2},e,f,{v\over2})=exp({-{1\over2}}\ad_v)\BCH(e,f)\\
&=e+f+{1\over48}(E^2f+F^2e)+\cdots
\end{align*}
where the terms with one and three Lie brackets vanish. This leads to the completed symmetric model of the bi-gon in which the first few orders of $\D{g}$ are
\begin{align*}
\D_0g&=e+f\\
\D_1g&=-{1\over2}[a+b,g]\\
\D_2g&={1\over48}(E^2f+F^2e)+{1\over16}[(E-F)(a-b),g]\\
\D_3g&=0
\end{align*}
This is distinct from the model given in [2]. Here $\D_0g$ is given by definition as the geometric boundary. Also $\D_1g$ is forced. The second order term is not unique, and the choice in [2] was different, there $\D_2g=\frac{1}{24}((F-E)G+G(F-E))(b-a)$. As predicted in [2], $\D_3g=0$ as for all odd orders $>1$.

\end{document}